\documentclass[10pt]{amsart}
\usepackage[initials]{amsrefs}
\usepackage{geometry,graphicx,amscd,color,amsmath,amsfonts,amssymb}
\numberwithin{equation}{section}
\begin{document}
\title[Lower bounds for the constants in the Bohnenblust-Hille inequality]{Lower bounds for the constants in the Bohnenblust-Hille inequality: the case of real scalars}
\author[D. Diniz \and G. A. Mu\~{n}oz-Fern\'{a}ndez \and D. Pellegrino \and J. B. Seoane-Sep\'{u}lveda]{D. Diniz, G. A. Mu\~{n}oz-Fern\'{a}ndez
\textsuperscript{*} \and D. Pellegrino \and J. B. Seoane-Sep\'{u}lveda\textsuperscript{*}}
\address{Unidade Academica de Matem\'{a}tica e Estat\'{i}stica,\newline\indent Universidade Federal de Campina Grande, \newline\indent Caixa Postal 10044, \newline\indent Campina Grande, 58429-970, Brazil.}
\email{diogo@dme.ufcg.edu.br}
\address{Departamento de An\'{a}lisis Matem\'{a}tico,\newline\indent Facultad de Ciencias Matem\'{a}ticas, \newline\indent Plaza de Ciencias 3, \newline\indent Universidad Complutense de Madrid,\newline\indent Madrid, 28040, Spain.}
\email{gustavo$\_$fernandez@mat.ucm.es}
\address{Departamento de Matem\'{a}tica, \newline\indent Universidade Federal da Para\'{\i}ba, \newline\indent 58.051-900 - Jo\~{a}o Pessoa, Brazil.} \email{pellegrino@pq.cnpq.br and dmpellegrino@gmail.com}
\address{Departamento de An\'{a}lisis Matem\'{a}tico,\newline\indent Facultad de Ciencias Matem\'{a}ticas, \newline\indent Plaza de Ciencias 3, \newline\indent Universidad Complutense de Madrid,\newline\indent Madrid, 28040, Spain.}
\email{jseoane@mat.ucm.es}
\subjclass[2010]{46G25, 47H60.}
\thanks{\textsuperscript{*} The second and fourth authors were supported by the Spanish Ministry of Science and Innovation (grant MTM2009-07848).}
\keywords{Absolutely summing operators, Bohnenblust-Hille Theorem.}
\begin{abstract}
The Bohnenblust-Hille inequality was obtained in 1931 and (in the case of real scalars) asserts that for every positive integer $N$ and every $m$-linear mapping $T:\ell_{\infty}^{N}\times\cdots\times\ell_{\infty}^{N}\rightarrow
\mathbb{R}$ one has
\begin{equation*}
\left( \sum\limits_{i_{1},...,i_{m}=1}^{N}\left\vert
T(e_{i_{^{1}}},...,e_{i_{m}})\right\vert ^{\frac{2m}{m+1}}\right) ^{\frac{m+1}{2m}}\leq C_{m}\left\Vert T\right\Vert ,
\end{equation*}
for some positive constant $C_{m}$. Since then, several authors obtained upper estimates for the values of $C_{m}$. However, the novelty presented in this short note is that we provide lower (and non-trivial) bounds for $C_{m}$.
\end{abstract}
\maketitle


\section{Introduction}

The Bohnenblust-Hille inequality (see \cite{bh}), for real scalars asserts
that for every positive integer $N$ and every $m$-linear mapping
$T:\ell_{\infty}^{N}\times\cdots\times\ell_{\infty}^{N}\rightarrow\mathbb{R}$
there is a positive constant $C_{m}$ such that
\begin{equation}
\left(  \sum\limits_{i_{1},...,i_{m}=1}^{N}\left\vert T(e_{i_{^{1}}%
},...,e_{i_{m}})\right\vert ^{\frac{2m}{m+1}}\right)  ^{\frac{m+1}{2m}}\leq
C_{m}\left\Vert T\right\Vert .\label{juui}%
\end{equation}
The Bohnenblust–Hille inequality has important applications in various fields of analysis. For details and references we mention \cite{defant}. When $m=2$ it is interesting to note that the Bohnenblust-Hille inequality is precisely the well-known Littlewood's $4/3$ inequality \cite{Litt}.

Since the 1930's many authors have obtained estimates for upper bounds of
$C_{m}$ in the case of real and complex scalars (see, e.g., \cite{bh, Davie,
defant2, defant, muno, Ka, Pell}). The constants of the polynomial version of the
Bohnenblust-Hille inequality (complex case) were recently investigated in \cite{dee}. Until
now, the more accurate upper bounds for $C_{m}$ (real case) were given in
\cite{Pell}:

\begin{center}
\begin{tabular}[c]{c|c}
$m$ & $\text{upper bounds for } C_{m}$\\
\hline
2 & $\sqrt{2} \approx1.414$\\
3 & $2^{\frac{20}{24}}\approx1.782$\\
4 & $2^{\frac{32}{32}}=2$\\
5 & $2^{\frac{48}{40}} \approx2.298$\\
6 & $2^{\frac{64}{48}} \approx2.520$\\
7 & $2^{\frac{84}{56}} \approx2.828$\\
8 & $2^{\frac{104}{64}} \approx3.084 $\\
9 & $2^{\frac{128}{72}} \approx3.429 $\\
10 & $2^{\frac{152}{80}} \approx3.732$
\end{tabular}
\end{center}
Also, it has very recently been proved (\cite{arxivpaper}) that the sequence of constants $(C_m)_m$ has the best possible asymptotic behavior, that is
\begin{equation*}
\displaystyle\lim_{m\rightarrow \infty }\frac{C_{m}}{C_{m-1}}=1.
\end{equation*}

However, and to the best of our knowledge, there is absolutely no work
presenting estimates for (non-trivial) lower bounds for the constants $C_{m}$.
In this short note we obtain lower bounds for $C_{m}$ which we believe are
(specially) interesting for the cases $m=2,3,4,5$.

In the following $e_k$ denotes the $k$th canonical vector in ${\mathbb R}^N$.

\section{The case $m=2$}\label{case2}

Let $T_{2}:\ell_\infty^2\times \ell_\infty^2\rightarrow\mathbb{R}$ be defined as
\[
T_{2}(x,y)=x_{1}y_{1}+x_{1}y_{2}+x_{2}y_{1}-x_{2}y_{2}.
\]
Note that $\left\Vert T_{2}\right\Vert =2.$ In fact,
\begin{align*}
|T_{2}(x,y)|  & =|x_{1}(y_{1}+y_{2})+x_{2}(y_{1}-y_{2})|\\
& \leq\Vert x\Vert(|y_{1}+y_{2}|+|y_{1}-y_{2}|)\\
& =2\Vert x\Vert\max\{|y_{1}|,|y_{2}|\}\\
& =2\Vert x\Vert\Vert y\Vert
\end{align*}
and since $T_{2}(e_{1}+e_{2},e_{1}+e_{2})=2$ it follows that $\left\Vert
T_{2}\right\Vert =2.$

Now the inequality
\[
\left(  \sum\limits_{i_{1},i_{2}=1}^{2}\left\vert T_{2}(e_{i_{^{1}}%
},e_{i_{2}})\right\vert ^{\frac{4}{3}}\right)  ^{\frac{3}{4}}\leq
C_{2}\left\Vert T_{2}\right\Vert
\]
can be re-written as
\[
4^{3/4}\leq2C_{2}%
\]
which gives
\[
C_{2}\geq2^{1/2}.
\]
Since it is well-known that $C_{2}\leq2^{1/2}$, we conclude that
$C_{2}=2^{1/2}$, but this result seems to be already known.

\section{The case $m=3$}\label{case3}

Now, let $T_{3}:\ell_\infty^4\times \ell_\infty^4\times \ell_\infty^4\rightarrow\mathbb{R}$ given by
\begin{align*}
T_{3}(x,y,z)   =&(z_{1}+z_{2})\left(  x_{1}y_{1}+x_{1}y_{2}+x_{2}y_{1}%
-x_{2}y_{2}\right) +(z_{1}-z_{2})\left(  x_{3}y_{3}+x_{3}y_{4}+x_{4}y_{3}-x_{4}y_{4}\right)  .
\end{align*}
We have
    \begin{align*}
    |T_3(x,y,z)| & =  |(z_{1}+z_{2})\left(  x_{1}y_{1}+x_{1}y_{2}+x_{2}y_{1}%
    -x_{2}y_{2}\right) +(z_{1}-z_{2})\left(  x_{3}y_{3}+x_{3}y_{4}+x_{4}y_{3}-x_{4}y_{4}\right)|\\
    & \leq  |z_{1}+z_{2}|\left(  |x_{1}||y_{1}+y_{2}|+|x_{2}||y_{1}
    -y_{2}|\right) +|z_{1}-z_{2}|\left(|x_{3}||y_{3}+y_{4}|+|x_{4}||y_{3}-y_{4}|\right)\\
    & \leq  \|x\|\Big\{|z_{1}+z_{2}|\left(|y_{1}+y_{2}|+|y_{1}
    -y_{2}|\right)
    +|z_{1}-z_{2}|\left(|y_{3}+y_{4}|+|y_{3}-y_{4}|\right)\Big\}\\
    & =  2\|x\|\Big\{|z_{1}+z_{2}|\max\{|y_1|,|y_2|\} +|z_{1}-z_{2}||\max\{|y_3|,|y_4|\}\Big\}\\
    & \leq  2\|x\| \|y\| \left(|z_{1}+z_{2}|+|z_{1}-z_{2}|\right)\\
    & =  4\|x\|\|y\|\max\{|z_1|,|z_2|\}\\
    &\leq 4\|x\|\|y\|\|z\|.
    \end{align*}
Since $T_3(e_1+e_2+e_3,e_1+e_2+e_3,e_1+e_2+e_3)=4$, then $\|T_3\|=4$. Also
\[
\left(  \sum\limits_{i_{1},i_2,i_{3}=1}^{4}\left\vert T_{3}(e_{i_{^{1}}%
},e_{i_2},e_{i_{3}})\right\vert ^{\frac{6}{4}}\right)  ^{\frac{4}{6}}\leq
C_{3}\left\Vert T_3\right\Vert
\]
becomes
\[
16^{2/3}\leq4C_{3}%
\]
which gives
\[
C_{3}\geq2^{2/3} \approx1.587.
\]

\section{The case $m=4$}

In this case, let us consider $T_{4}:\ell_\infty^8\times \ell_\infty^8\times \ell_\infty^8\times
\ell_\infty^8\rightarrow\mathbb{R}$ given by
\begin{align*}
&  T_{4}(x,y,z,w) =\\
&  =\left(  w_{1}+w_{2}\right)  \left(
\begin{array}
[c]{c}%
(z_{1}+z_{2})\left(  x_{1}y_{1}+x_{1}y_{2}+x_{2}y_{1}-x_{2}y_{2}\right) \\
+(z_{1}-z_{2})\left(  x_{3}y_{3}+x_{3}y_{4}+x_{4}y_{3}-x_{4}y_{4}\right)
\end{array}
\right) \\
&  +\left(  w_{1}-w_{2}\right)  \left(
\begin{array}
[c]{c}%
(z_{3}+z_{4})\left(  x_{5}y_{5}+x_{5}y_{6}+x_{6}y_{5}-x_{6}y_{6}\right) \\
+(z_{3}-z_{4})\left(  x_{7}y_{7}+x_{7}y_{8}+x_{8}y_{7}-x_{8}y_{8}\right).
\end{array}
\right)
\end{align*}
As in Sections \ref{case2} and \ref{case3} we see that $\left\Vert T_{4}\right\Vert =8$ and from (\ref{juui}) we obtain
\[
64^{5/8}\leq8 C_{4}.%
\]
Hence
\[
C_{4}\geq2^{3/4}\approx1.681.
\]

%

\section{The general case}

From the previous results it is not difficult to prove that in general
\[
C_{m}\geq2^{\frac{m-1}{m}}
\]
for every $m\geq 2$.
Indeed, let us define the $m$-linear forms $T_m:\ell_\infty^{2^{m-1}}\times\stackrel{m}{\ldots}\times\ell_\infty^{2^{m-1}}\rightarrow{\mathbb R}$ by induction as
    \begin{align*}
    T_2(x_1,x_2)=&x^1_1x^1_2+x^1_1x^2_2+x^2_1x^1_2-x^2_1x^2_2\\
    T_m(x_1,\ldots,x_m)=&(x_m^1+x_m^2)T_{m-1}(x_1,\ldots,x_{m-1})\\
    &+(x_m^1-x_m^2)T_{m-1}(B^{2^{m-2}}(x_1),B^{2^{m-2}}(x_2),B^{2^{m-3}}(x_3)\ldots,B^2(x_{m-1})),
    \end{align*}
where $x_k=(x_k^n)_n\in \ell_\infty^{2^{m-1}}$ for $1\leq k\leq m$, $1\leq n\leq 2^{m-1}$ and $B$ is the backward shift operator in $\ell_\infty^{2^{m-1}}$. It has been proved in Section \ref{case2}
that $\|T_2\|=2$. If we assume that $\|T_{m-1}\|=2^{m-2}$, then
    \begin{align*}
    |T_m(x_1,\ldots,x_m)\leq & |x_m^1+x_m^2||T_{m-1}(x_1,\ldots,x_{m-1})|\\
    &+|x_m^1-x_m^2||T_{m-1}(B^{2^{m-2}}(x_1),B^{2^{m-2}}(x_2),B^{2^{m-3}}(x_3)\ldots,B^2(x_{m-1}))|\\
    \leq & 2^{m-2}[|x_m^1+x_m^2|\|x_1\|\cdots\|x_{m-1}\|\\
    &+|x_m^1-x_m^2|\|B^{2^{m-2}}(x_1)\|\|B^{2^{m-2}}(x_2)\|\|B^{2^{m-1}}(x_3)\|\cdots\|B^{2}(x_{m-1})\|]\\
    \leq & 2^{m-2}[|x_m^1+x_m^2|+|x_m^1-x_m^2|]\|x_1\|\cdots\|x_{m-1}\|\\
    =&2^{m-1}\|x_1\|\cdots\|x_{m-1}\|\max\{|x_m^1|,|x_m^2|\}\\
    \leq & 2^{m-1}\|x_1\|\cdots\|x_m\|.
    \end{align*}
This induction argument shows that $\|T_m\|\leq 2^{m-1}$ for all $m\in{\mathbb N}$. Using a similar induction argument it is easy to prove that
$T_m(x_1,\ldots,x_m)=2^{m-1}$ for $x_1,\ldots,x_m$ such that
$\|x_1\|=\ldots=\|x_m\|=1$ and $x^k_j=1$ with $1\leq j,k\leq m$, which proves that $\|T_m\|=2^{m-1}$ for all $m\in{\mathbb N}$.

On the other hand from \eqref{juui} we have
    $$
    \left(  \sum\limits_{i_{1},...,i_{m}=1}^{2^{m-1}}\left\vert T_{m}(e_{i_{^{1}}%
    },...,e_{i_{m}})\right\vert ^{\frac{2m}{m+1}}\right)  ^{\frac{m+1}{2m}}\leq
    C_{m}\left\Vert T_m\right\Vert=2^{m-1}C_m.
    $$
To finish we shall prove that  $|T_{m}(e_{i_{^{1}}%
    },...,e_{i_{m}})|$ is either $0$ or $1$ and that $|T_{m}(e_{i_{^{1}}%
    },...,e_{i_{m}})|=1$ for exactly $4^{m-1}$ choices of the vectors $e_{i_1},\ldots,e_{i_m}$. Working again by induction, the reader can easily check that the latter is true for $m=2$ (see Section \ref{case2}). If we assume that the result is true for $m-1$ and $e_{i_1},\ldots,e_{i_{m-1}}$ is one of the $4^{m-2}$ choices of the unit vectors such that $|T_{m-1}(e_{i_1},\ldots,e_{i_{m-1}})|=1$, then using the definition of $T_m$
        \begin{align*}
            |T_m(e_{i_1},\ldots,e_{i_{m-1}},e_k)|&=|T_{m-1}(e_{i_1},\ldots,e_{i_{m-1}})|=1,\\
            |T_m(e_{i_1+2^{m-2}},e_{i_2+2^{m-2}},e_{i_3+2^{m-3}}\ldots,e_{i_{m-1}+2},e_k)|&=|T_{m-1}(e_{i_1},\ldots,e_{i_{m-1}})|=1,
        \end{align*}
for $k=1,2$. Hence we have found $4\cdot 4^{m-2}=4^{m-1}$ choices of unit vector at which $|T_m|$ takes the value 1. A simple inspection of the problem shows that $|T_m|$ vanishes at any other choice of canonical vectors.

\section{Final remarks}

Notice that our estimate $2^\frac{m-1}{m}$ seems inaccurate as $m\rightarrow\infty$ since it is a common
feeling that the optimal values for the constants $C_{m}$ should tend to
infinity as $m\rightarrow\infty$. However, and as a matter of fact, we must
say that this common feeling seems to be just supported by the estimates of
the upper bounds for $C_{m}$ obtained throughout the last decades, but there
seems to be not any particular result supporting this ``fact''. In any case (and at
least for $m=2,3,4,5$) our estimates are clearly interesting. Summarizing:

\begin{center}%
\begin{tabular}
[c]{rcl}
& $C_{2}$ & $=\sqrt{2}$\\
$1.587 \leq$ & $C_{3}$ & $\leq1.782$\\
$1.681 \leq$ & $C_{4}$ & $\leq2$\\
$1.741 \leq$ & $C_{5}$ & $\leq2.298$.
\end{tabular}

\end{center}

We also conclude that $C_{3}>C_{2}$, which seems to be not known until now.


\end{document}